\documentclass[12pt]{article}

\usepackage{amsmath,amssymb}
\usepackage[applemac]{inputenc}
\usepackage{amsmath,amssymb,fullpage}
\usepackage{showlabels}
\usepackage{color}

\newcommand{\dproof}{\noindent {Proof.} \quad}
\newcommand{\fproof}{\hfill $\square$ \bigskip}

\newtheorem{definition}{Definition}[section]

\newtheorem{theorem}[definition]{Theorem}
\newtheorem{problem}[definition]{Problem}
\newtheorem{remark}[definition]{ \it Remark}
\newtheorem{coro}[definition]{Corollary}

\numberwithin{equation}{section}

\def\1B{\text{1\!\!I}}

\begin{document}
\date{1 June 2018 }
\title{A white noise approach to optimal insider control of systems with delay}

\author{
Olfa Draouil$^{1,2}$ and Bernt \O ksendal$^{1,2}$}

\footnotetext[1]{Department of Mathematics, University of Oslo, P.O. Box 1053 Blindern, N--0316 Oslo, Norway.\\
Email: {\tt olfad@math.uio.no\\ \tt oksendal@math.uio.no}}

\footnotetext[2]{This research was carried out with support of the Norwegian Research Council, within the research project Challenges in Stochastic Control, Information and Applications (STOCONINF), project number 250768/F20.}

\maketitle
\paragraph{MSC(2010):} 60H05; 60H07; 60H40; 60G57; 91B70; 93E20.

\paragraph{Keywords:} Stochastic delay equation, optimal insider control, Hida-Malliavin derivative, Donsker delta functional, white noise theory, stochastic maximum principles, time-advanced BSDE, optimal insider portfolio in a financial market with delay.

\begin{abstract} We use a white noise approach to study the problem of optimal inside control of a stochastic delay equation driven by a Brownian motion $B$ and a Poisson random measure $N$. In particular, we use Hida-Malliavin calculus and the Donsker delta functional to study the problem.
\vskip 0.2cm
We establish a sufficient and a necessary maximum principle for the optimal control when the trader from the beginning has inside information about the future value of some random variable related to the system.These results are applied to the problem of optimal inside harvesting control in a population modelled by a stochastic delay equation.
\vskip 0.2cm
Next, we apply a direct white noise method to find the optimal insider portfolio in a financial market where the risky asset price is given by a stochastic delay equation. A classical result of Pikovski and Karatzas shows that when the inside information is $B(T)$, where $T$ is the terminal time of the trading period, then the market is not viable.
Our results show that with this inside information the market is not viable even if there is delay in the equations.

\end{abstract}

\section{Introduction}
Delay stochastic differential equations are a type of stochastic differential equation in which the derivative of the solution process at any given time depends not only on its value at the present time, but also on its values at previous times.
They are also called systems with aftereffect.

Many real life systems in for example engineering, biology, finance and communication networks include delay phenomena in their dynamics.
Therefore delay stochastic differential equations have many applications, and the interest of  delay stochastic differential equations continues to grow in all areas of science and in particular in control engineering.

In our paper, in addition to the delay, we suppose that we have an additional information related to the future value of the system. In other words, we consider insider control of delay stochastic differential equations.
Insider control problems are treated in many papers. We mention \cite{AIS,A,AZ,JM,JY}.
Our paper differs from the above papers in two ways:
\begin{itemize}
\item
 We use a different insider control approach based on white noise theory. Specifically, we apply Hida-Malliavin calculus and the definition of the Donsker delta functional. This allows us to transform the insider delay control problem into a non insider parameterized control problem. We prove two maximum principle theorems. Then we apply the results to an inside optimal harvesting problem in a population modelled by a delay equation.


 \item
We apply a direct white noise method to study optimal inside control in the context of delay. In particular, an interesting question which, to the best of our knowledge, has not been studied before, is how inside information combines with the delay in the dynamics when it comes to performance, for example in finance. Specifically, we find an explicit expression for the optimal insider portfolio in a delay market and we prove that, perhaps surprisingly, even in the presence of the delay, the market is not viable when the inside information is given by $B(T)$, where $T$ is the terminal value of the trading interval.   See Subsection 7.1.
 \end{itemize}
We now explain this in more detail.

In this paper we consider an insider's optimal control problem for a stochastic process $X(t)=X(t,Z)=X^{u}(t,Z)$ defined as the solution of a stochastic differential delay equation  of the form

\begin{equation}\label{eq1.1bis}
\begin{cases}
dX(t)=dX(t,Z)=b(t,X(t,Z),Y(t,Z),u(t,Z),Z)dt+\sigma(t,X(t,Z),Y(t,Z),u(t,Z),Z)dB(t)\\
+\int_{\mathbb{R}}\gamma(t,X(t,Z),Y(t,Z),u(t,Z),Z,\zeta)\tilde{N}(dt,d\zeta), \quad 0\leq t\leq T\\
X(t)=\xi(t), \quad -\delta\leq t\leq 0,
\end{cases}
\end{equation}
where
\begin{equation}
Y(t,Z)=X(t-\delta,Z),
\end{equation}
$\delta > 0$ being a fixed constant (the delay) and $\xi$ is a deterministic function.\\
Here $B(t)$ and $\tilde{N}(dt,d\zeta)$ is a Brownian motion and an independent compensated Poisson random measure, respectively, jointly defined on a filtered probability space $(\Omega, \mathbb{F}=\{ \mathcal{F}_t \}_{t \in [0,T]},\mathbf{P})$ satisfying the usual conditions where $\Omega=S'(\mathbb{R})$ the dual of Schwartz space and $\mathbf{P}$  is the Gaussian measure on $S'(\mathbb{R})$. $T>0$ is a given constant. We refer to \cite{DOP} for more information about white noise theory and \cite{OS1} for  stochastic calculus for It\^{o}-L\' evy processes.

The process $u(t,Z)=u(t,x,z)_{z=Z}$ is our insider control process, where $Z$ is a given $\mathcal{F}_{T_0}$-measurable random variable for some $T_0>T$ , representing the inside information available to the controller.

We assume that the inside information is of \emph{initial enlargement} type. Specifically, we assume  that the inside filtration $\mathbb{H}$ has the form

\begin{equation}\label{eq1.1}
 \mathbb{H}= \{ \mathcal{H}_t\}_{0\leq t \leq T}, \text{ where } \mathcal{H}_t = \mathcal{F}_t \vee \sigma(Z)
\end{equation}
for all $t\in[0,T]$, where $Z$ is a given $\mathcal{F}_{T_0}$-measurable random variable, for some $T_0 > T$ (constant).
Here and in the following we use the right-continuous version of  $ \mathbb{H}$, i.e. we put
$\mathcal{H}_{t}= \mathcal{H}_{t^+}=\bigcap_{s>t}\mathcal{H}_s.$

We assume that the value at time $t$ of our insider control process $u(t)$ is allowed to depend on both $Z$ and $\mathcal{F}_t$. In other words, $u(.)$ is assumed to be $\mathbb{H}$-adapted,
such that $u(. , z)$ is $\mathbb{F}$-adapted for each $z \in \mathbb{R}$.

We also assume that the \emph{Donsker delta functional} of $Z$ exists. This assumption implies that the Jacod condition holds, and hence that $B(\cdot)$ and $N(\cdot,\cdot)$ are semimartingales with respect to $\mathbb{H}$ therefore equation \eqref{eq1.1bis} is well defined. We will explain this  with more details in the next Section.

Let $\mathbb{U}$ denote the set of admissible control values. We assume that the functions
\begin{align}
b(t,x,y,u,z)&=b(t,x,y,u,z,\omega):
[0,T]\times \mathbb{R} \times \mathbb{R} \times \mathbb{U}\times \mathbb{R} \times \Omega \mapsto \mathbb{R},\nonumber\\
\sigma(t,x,y,u,z)&=\sigma(t,x,y,u,z,\omega):
[0,T]\times \mathbb{R}\times \mathbb{R} \times \mathbb{U} \times \mathbb{R} \times \Omega \mapsto \mathbb{R},\nonumber\\
\gamma(t,x,y,u,z,\zeta)&=\gamma(t,x,y,u,z,\zeta,\omega):
[0,T]\times \mathbb{R} \times \mathbb{R} \times \mathbb{U} \times \mathbb{R} \times \mathbb{R} \times \Omega \mapsto \mathbb{R},\nonumber
\end{align}
are given  $C^1$ functions with respect to $x$, $y$ and $u$ and adapted processes in $(t,\omega)$ for each given $x,y,u,z,\zeta$.
Let $\mathcal{A}$ be a given family of admissible $\mathbb{H}-$adapted controls $u$.
The \emph{performance functional} $J(u)$ of a control process $u \in \mathcal{A}$ is defined by
\begin{equation}\label{eq1.4}
J(u)= \mathbb{E}[\int_0^T f(t,X(t,Z),u(t,Z),Z))dt +g(X(T,Z),Z)],
\end{equation}
where \begin{align}
&f(t,x,u,z): [0,T]  \times\mathbb{R}\times \mathbb{U} \times\mathbb{R} \mapsto \mathbb{R}\nonumber\\
&g(x,z):\mathbb{R}\times\mathbb{R} \mapsto \mathbb{R}
\end{align}
are given functions, $C^1$ with respect to $x$ and $u$. The functions $f$ and $g$ are called the \emph{profit rate} and \emph{terminal payoff}, respectively. For completeness of the presentation we allow these functions to depend explicitly on the future value $Z$ also, although this would not be the typical case in applications. But it could be that $f$ and $g$ are influenced by the future value $Z$ directly through the action of an insider, in addition to being influenced indirectly through the control process $u$ and the corresponding state process $X$.\\
 The problem we consider is the following:

\begin{problem}
Find $u^{\star} \in\mathcal{A}$ such that
\begin{equation}\label{eq1.5}
    \sup_{u\in\mathcal{A}}J(u)=J(u^{\star}).
\end{equation}
\end{problem}

This paper is organized as follow:
\begin{itemize}
\item In Section 2 we recall the definition of the Donsker delta functional and its basic properties.
\item In Section 3 we use the Donsker delta functional to transform the insider control problem to a parametrized non-insider problem.
\item In Section 4 and 5  we prove a sufficient and necessary maximum principle theorems.
\item In Section 6 we use the results of the previous sections to study an inside optimal harvesting problem in a population modelled by delay equations.
\item In Section 7 we apply a direct white noise method to study the problem of optimal insider portfolio in a financial market with delay. This example is not directly related to the previous sections since we do not use the maximum principle theorems in order to solve it. Then we study the viability of this market with delay.
\end{itemize}

\section{The Donsker delta functional}
To study the Problem 1.1 we adapt a white noise approach, in particular we will use the definition of the Donsker delta funtional which takes value in the stochastic Hida distribution space $(\mathcal{S})^*$ and the Hida-Malliavin derivatives  $D_t$ and $D_{t,z}$ with respect to the Brownian motion $B(.)$ and the Poisson random measure $\tilde{N}(.,.)$, respectively. We refer to  \cite{DOP} for more information about $(\mathcal{S})^*$, white noise theory and the Hida-Malliavin derivatives.
Let us now give the definition and basic properties of the Donsker delta functional:

\begin{definition}
Let $Z:\Omega\rightarrow\mathbb{R}$ be a random variable which also belongs to $(\mathcal{S})^{\ast}$. Then a continuous functional
\begin{equation}\label{donsker}
    \delta_{Z}(.): \mathbb{R}\rightarrow (\mathcal{S})^{\ast}
\end{equation}
is called a Donsker delta functional of $Z$ if it has the property that
\begin{equation}\label{donsker property }
    \int_{\mathbb{R}}g(z)\delta_{Z}(z)dz= g(Z) \quad a.s.
\end{equation}
for all (measurable) $g : \mathbb{R} \rightarrow \mathbb{R}$ such that the integral converges.
\end{definition}
\begin{remark}
 If $Z\in L^2(S'(\mathbb{R}),\mathbf{P})$ then it belongs to $(\mathcal{S})^*$.
\end{remark}
\subsection{Donsker delta functional properties}
Define the \emph{regular conditional distribution} with respect to $\mathcal{F}_t$ of a given real random variable $Z$, denoted by $Q_t(dz)=Q_t(\omega,dz)$, by the following properties:
\begin{itemize}
\item
For any Borel set $\Lambda \subseteq \mathbb{R}, Q_t(\cdot, \Lambda)$ is a version of $\mathbb{E}[\mathbf{1}_{Z \in \Lambda} | \mathcal{F}_t]$.
\item
For each fixed $\omega, Q_t(\omega,dz)$ is a probability measure on the Borel subsets of $\mathbb{R}$.
\end{itemize}
Our aim is to show that if Z has a Donsker delta functional then Jacod condition holds hence $B(.)$ and  and $N(.,.)$ are semimartingales with respect
to $\mathbb{H}$. Therefore the SDDE \eqref{eq1.1bis} is well defined.

First let us recall the Jacod condition:\cite{AJ, Jac}\\
If there exists a $\sigma$-finite positive measure $\eta$ such that the regular conditional distributions
of the random variable $Y$ given $\mathcal{F}_t, t\in[0, T]$ are absolutely continuous with respect to
$\eta$, Jacod proves that every $(P, \mathbb{F})$-martingale remains a $(P,\mathbb{H})$-semimartingale on the
interval $[0, T]$.

Now assume that the Donsker delta functional of $Z$ exists. Equation \eqref{donsker property } holds for all
bounded $g$. Taking the conditional expectation of \eqref{donsker property } with respect to $\mathcal{F}_t$ we get
\begin{align}
\int_{\mathbb{R}}g(z)E[\delta_Z (z)|\mathcal{F}_t]dz &= E[g(Z )|\mathcal{F}_t]\nonumber\\
&=\int_{\mathbb{R}} g(z)Q_t(\omega, dz)
\end{align}
Since this holds for all $g$, we conclude that we have the following identity of measures:
\begin{equation}
Q_t(\omega, dz) = E[\delta_Z (z)|\mathcal{F}_t]dz.
\end{equation}
Hence $Q_t(\omega, dz)$ is absolutely continuous with respect to the Lebesgue measure. In particular,
the Jacod condition holds then $B(.)$ and $N(.,.)$ are semimartingales with respect
to $\mathbb{H}$. Therefore the SDDE \eqref{eq1.1bis} is well defined.\\
In the next subsection we will give some examples of Donsker delta functional.
\subsection{The Donsker delta functional for a class of It\^{o} - L\'{e}vy processes}

Consider the special case when $Z$ is a first order chaos random variable of the form
\begin{equation}\label{eq2.5}
    Z = Z (T_0); \text{ where } Z (t) =\int_0^t\beta(s)dB(s)+\int_0^t\int_{\mathbb{R}}\psi(s,\zeta)\tilde{N}(ds,d\zeta), \mbox{ for } t\in [0,T_0]
\end{equation}
for some deterministic functions $\beta \neq 0, \psi$ such that
\begin{equation}\label{}
    \int_0^{T_0} \{ \beta^2(t)+\int_{\mathbb{R}}\psi^2(t,\zeta)\nu(d\zeta)\} dt<\infty \text{ a.s. }
\end{equation}
and for every $\epsilon >0$ there exists $\rho > 0$ such that
\begin{equation}\label{eq8.4}
\int_{\mathbb{R} \setminus (-\epsilon,\epsilon)} e^{\rho  \zeta} \nu(d\zeta) < \infty.\nonumber\\
\end{equation}

This condition implies that the polynomials are dense in $L^2(\mu)$, where $d\mu(\zeta)=\zeta^2 d\nu(\zeta)$. It also guarantees that the measure $\nu$ integrates all polynomials of degree $\geq 2$.\\
In this case it is well known (see e.g. \cite{MOP}, \cite{DiO1}, Theorem 3.5, and \cite{DOP},\cite{DiO2}) that the Donsker delta functional exists in $(\mathcal{S})^{\ast}$ and is given
by
\begin{eqnarray}\label{eq2.7}
   \delta_Z(z)&=&\frac{1}{2\pi}\int_{\mathbb{R}}\exp^{\diamond}\big[ \int_0^{T_0}\int_{\mathbb{R}}(e^{ix\psi(s,\zeta)}-1)\tilde{N}(ds,d\zeta)+ \int_0^{T_0}ix\beta(s)dB(s)  \nonumber\\
   &+&  \int_0^{T_0}\{\int_{\mathbb{R}}(e^{ix\psi(s,\zeta)}-1-ix\psi(s,\zeta))\nu(d\zeta)-\frac{1}{2}x^2\beta^2(s)\}ds-ixz\big]dx,
\end{eqnarray}
where $\exp^{\diamond}$ denotes the Wick exponential.
Moreover, we have for $t< T_0$
\begin{align}
&\mathbb{E}[\delta_Z(z)|\mathcal{F}_t]\nonumber\\
=&  \frac{1}{2\pi}\int_{\mathbb{R}}\exp\big[\int_0^t\int_{\mathbb{R}}ix\psi(s,\zeta)\tilde{N}(ds,d\zeta) +\int_0^t ix\beta(s)dB(s)\\
&+\int_t^{T_0}\int_{\mathbb{R}}(e^{ix\psi(s,\zeta)}-1-ix\psi(s,\zeta))\nu(d\zeta)ds-\int_t^{T_0}\frac{1}{2}x^2\beta^2(s)ds-ixz\big]dx.
\end{align}

If $D_t$ and $D_{t,\zeta}$ denotes the \emph{Hida-Malliavin derivative} at $t$ and $t,\zeta$  with respect to $B$ and $\tilde{N}$, respectively, we have
\begin{align}
&\mathbb{E}[D_t\delta_Z(z)|\mathcal{F}_t]=\nonumber\\
&\frac{1}{2\pi}\int_{\mathbb{R}}\exp\big[\int_0^t\int_{\mathbb{R}}ix\psi(s,\zeta)\tilde{N}(ds,d\zeta) +\int_0^t ix\beta(s)dB(s)\nonumber\\
&+\int_t^{T_0}\int_{\mathbb{R}}(e^{ix\psi(s,\zeta)}-1-ix\psi(s,\zeta))\nu(d\zeta)ds-\int_t^{T_0}\frac{1}{2}x^2\beta^2(s)ds-ixz\big]ix\beta(t)dx
\end{align}
   and
\begin{align}
&\mathbb{E}[D_{t,z}\delta_Z(z)|\mathcal{F}_t]=\nonumber\\
& \frac{1}{2\pi}\int_{\mathbb{R}}\exp\big[\int_0^t\int_{\mathbb{R}}ix\psi(s,\zeta)\tilde{N}(ds,d\zeta) +\int_0^t ix\beta(s)dB(s)\nonumber\\
&+\int_t^{T_0}\int_{\mathbb{R}}(e^{ix\psi(s,\zeta)}-1-ix\psi(s,\zeta))\nu(d\zeta)ds-\int_t^{T_0}\frac{1}{2}x^2\beta^2(s)ds-ixz\big](e^{ix\psi(t,z)}-1)dx.
\end{align}

\subsection{The Donsker delta functional for a Gaussian process}

 Consider the special case when $Z$ is a Gaussian random variable of the form
\begin{equation}\label{eq5.47}
    Z = Z (T_0); \text{ where } Z (t) =\int_0^t\beta(s)dB(s), \mbox{ for } t\in [0,T_0],
\end{equation}
for some deterministic function $\beta\in \mathbf{L}^2[0,T_0]$ with
\begin{equation}\label{}
    \|\beta\|^2_{[t,T]} :=\int_t^T\beta(s)^2ds>0 \mbox{ for all } t\in[0,T].
\end{equation}
In this case it is well known that the Donsker delta functional is given
by
\begin{equation}\label{}
    \delta_{Z}(z)=(2\pi v)^{-\frac{1}{2}}\exp^{\diamond}[-\frac{(Z-z)^{\diamond2}}{2v}],
\end{equation}
where we have put $v :=\|\beta\|^2_{[0,T_0]}$. See e.g. \cite{AaOU}, Proposition $3.2$.
Using the Wick rule when taking conditional expectation, using the martingale property of
the process $Z (t)$ and applying Lemma $3.7$ in \cite{AaOU} we get
\begin{eqnarray}\label{eq5.50}
   \mathbb{E}[\delta_Z(z)|\mathcal{F}_t]&=&(2\pi v)^{-\frac{1}{2}}\exp^{\diamond}[-\mathbb{E}[\frac{(Z(T_0)-z)^{\diamond 2}}{2v}|\mathcal{F}_t]] \nonumber \\
   &=& (2\pi \|\beta\|^2_{[0,T_0]})^{-\frac{1}{2}}\exp^{\diamond}[- \frac{(Z(t)-z)^{\diamond 2}}{2\|\beta\|^2_{[0,T_0]}}]\nonumber \\
   &=& (2\pi \|\beta\|^2_{[t,T_0]})^{-\frac{1}{2}} \exp[- \frac{(Z(t)-z)^2}{2\|\beta\|^2_{[t,T_0]}}].
\end{eqnarray}
Similarly, by the Wick chain rule and Lemma $3.8$ in \cite{AaOU} we get, for $t \in [0,T],$
\begin{eqnarray}\label{eq5.51}
  \mathbb{E}[D_t\delta_Z(z)|\mathcal{F}_t] &=&-\mathbb{E}[(2\pi v)^{-\frac{1}{2}}\exp^{\diamond}[- \frac{(Z(T_0)-z)^{\diamond 2}}{2v}]\diamond\frac{Z(T_0)-z}{v}\beta(t)|\mathcal{F}_t] \nonumber\\
   &=&-(2\pi v)^{-\frac{1}{2}} \exp^{\diamond}[- \frac{(Z(t)-z)^{\diamond 2}}{2v}]\diamond\frac{Z(t)-z}{v}\beta(t)\nonumber \\
   &=& -(2\pi \|\beta\|^2_{[t,T_0]})^{-\frac{1}{2}}\exp[- \frac{(Z(t)-z)^2}{2\|\beta\|^2_{[t,T_0]}}]\frac{Z(t)-z}{\|\beta\|^2_{[t,T_0]}}\beta(t).
\end{eqnarray}

For more information about the Donsker delta functional, Hida-Malliavin calculus and their properties, see \cite{DO1}.

From now on we assume that $Z$ is a given random variable which also belongs to $(\mathcal{S})^{\ast}$, with a Donsker delta functional $\delta_Z(z)\in(\mathcal{S})^{\ast}$
satisfying
\begin{equation}
\mathbb{E}[\delta_Z(z)|\mathcal{F}_T] \in \mathbf{L}^2(\mathcal{F}_T,P).
\end{equation}

\section{Transforming the insider control problem to a related parametrised non-insider problem}
Since $X(t)$ is $\mathbb{H}$-adapted, we get by using the definition of the Donsker delta functional $\delta_Z(z)$ of $Z$ that
\begin{equation}\label{eq1.6}
X(t)=X(t,Z)=X(t,z)_{z=Z}=\int_{\mathbb{R}}X(t,z)\delta_Z(z)dz,
\end{equation}
for some $z$-parametrised process $X(t,z)$ which is $\mathbb{F}$-adapted for each $z$.

Then, again by the definition of the Donsker delta functional we can write, for $0\leq t\leq T$
\begin{align}\label{eq1.7}
&X(t)= \xi(0) +\int_0^t b(s,X(s),Y(s),u(s,Z),Z)ds + \int_0^t \sigma(s,X(s),Y(s),u(s,Z),Z)dB(s)\nonumber\\
&+\int_0^t \int_{\mathbb{R}} \gamma(s,X(s),Y(s),u(s,Z),Z,\zeta)\tilde{N}(ds,d\zeta)\nonumber\\
&=\xi(0)+\int_0^tb(s,X(s,z),Y(s,z),u(s,z),z)_{z=Z}ds \nonumber\\
&+ \int_0^t \sigma(s,X(s,z),Y(s,z),u(s,z),z)_{z=Z}dB(s)\nonumber\\
&+\int_0^t \int_{\mathbb{R}} \gamma(s,X(s,z),Y(s,z),u(s,z),z,\zeta)_{z=Z}\tilde{N}(ds,d\zeta)\nonumber\\
&= \int_{\mathbb{R}} \xi(0)\delta_Z(z)dz+\int_0^t \int_{\mathbb{R}}b(s,X(s,z),Y(s,z),u(s,z),z)\delta_Z(z)dzds\nonumber\\
&+ \int_0^t \int_{\mathbb{R}}\sigma(s,X(s,z),Y(s,z),u(s,z),z)\delta_Z(z)dzdB(s)\nonumber\\
&+\int_0^t \int_{\mathbb{R}}\int_{\mathbb{R}} \gamma(s,X(s,z),Y(s,z),u(s,z),z,\zeta)\delta_Z(z)dz\tilde{N}(ds,d\zeta)\nonumber\\
&=\int_{\mathbb{R}} \{\xi(0)+\int_0^t b(s,X(s,z),Y(s,z),u(s,z),z)ds + \int_0^t \sigma(s,X(s,z),Y(s,z),u(s,z),z)dB(s)\nonumber\\
&+\int_0^t \int_{\mathbb{R}} \gamma(s,X(s,z),Y(s,z),u(s,z),z,\zeta)\tilde{N}(ds,d\zeta)\} \delta_Z(z)dz.
\end{align}

Comparing \eqref{eq1.6} and \eqref{eq1.7}
we see that  \eqref{eq1.6} holds if we for each $z$ choose $X(t,z)$ as the solution of the classical (but parameterized) SDDE
\begin{equation} \label{eq3.3}
\begin{cases}
dX(t,z) = b(t,X(t,z),Y(t,z),u(t,z),z)dt + \sigma(t,X(t,z),Y(t,z),u(t,z),z)dB(t)\\
    + \int_{\mathbb{R}} \gamma(t,X(t,z),Y(t,z),u(t,z),z,\zeta)\tilde{N}(dt,d\zeta); \quad t\in [0,T]\\
X(t,z)=\xi(t); \quad t \in [-\delta, 0].
\end{cases}
\end{equation}

For results about existence and uniqueness of solutions of SDDE  we refer to \cite{MS} and the references therein.

As before let $\mathcal{A}$ be the given family of admissible $\mathbb{H}-$adapted controls $u$.
Then in terms of $X(t,z)$ the performance functional $J(u)$ of a control process $u \in \mathcal{A}$ gets the form
\begin{align}\label{eq0.13}
J(u)&=\mathbb{E}[\int_0^T f(t,X(t,Z),u(t,Z),Z)dt +g(X(T,Z),Z)]\nonumber\\
&=\mathbb{E}[\int_0^T (\int_{\mathbb{R}}f(t,X(t,z),u(t,z),z)\mathbb{E}[\delta_Z(z)|\mathcal{F}_t]dz)dt +\int_{\mathbb{R}}g(X(T,z),z)\mathbb{E}[\delta_Z(z)|\mathcal{F}_T]dz]\nonumber\\
&= \int_{\mathbb{R}} j(u)(z) dz,
\end{align}
where
\begin{align}\label{eq1.5}
j(u)(z)&:=  \mathbb{E}[ \int_0^T f(t, X(t,z),u(t,z),z)\mathbb{E}[\delta_Z(z)|\mathcal{F}_t]dt \nonumber \\
   & + g(X(T,z),z)\mathbb{E}[\delta_Z(z)|\mathcal{F}_T]].
 \end{align}

 Thus we see that to maximise $J(u)$ it suffices to maximise $j(u)(z)$ for each value of the parameter $z \in \mathbb{R}$. Therefore Problem 1.1 is transformed into the problem
\begin{problem}
For each given $z \in \mathbb{R}$ find $u^{\star}=u^{\star}(t,z) \in\mathcal{A}$ such that
\begin{equation}\label{problem2}
    \sup_{u\in\mathcal{A}}j(u)(z)=j(u^{\star})(z).
\end{equation}
\end{problem}

\section{A sufficient-type maximum principle}
In this section we will establish a sufficient maximum principle for Problem 3.1.\\

Problem 3.1 is a stochastic control problem with a standard (parametrised) stochastic delay differential equation \eqref{eq3.3} for the state process $X(t,z)$, but with a non-standard performance functional given by \eqref{eq1.5}. We can solve this problem by a modified maximum principle approach, as follows:\\

Define the \emph{Hamiltonian}
$ H:[0,T]\times \mathbb{R}\times\mathbb{R}\times\mathbb{U}\times\mathbb{R}\times\mathbb{R}\times\mathbb{R}\times\mathcal{R} \times \Omega \rightarrow \mathbb{R}$ by
\begin{align}\label{eq4.1}
&H(t,x,y,u,z,p,q,r)=H(t,x,y,u,z,p,q,r,\omega)\nonumber\\
&=\mathbb{E}[\delta_Z(z)|\mathcal{F}_t] f(t,x,u,z)+b(t,x,y,u,z)p\nonumber\\
& + \sigma(t,x,y,u,z)q+\int_{\mathbb{R}}\gamma(t,x,y,u,z,\zeta)r(\zeta)\nu(d\zeta).
\end{align}

 $\mathcal{R}$ denotes the set of all functions $r(\cdot) : \mathbb{R}\rightarrow  \mathbb{R}$
such that the last integral above converges.
The quantities $p,q,r(\cdot)$ are called the \emph{adjoint variables}.
The \emph{adjoint processes} $p(t,z),q(t,z),r(t,z,\zeta)$ are defined as the solution of the $z$-parametrized advanced backward stochastic differential equation (ABSDE)

\begin{equation}\label{eq4.2a}
    \left\{
\begin{array}{l}
dp(t,z)=\mathbb{E}[\mu(t,z)|\mathcal{F}_t]dt+q(t,z) dB(t)+\int_{\mathbb{R}}r(t,z,\zeta)\tilde{N}(dt,d\zeta), \quad t \in[0,T]\nonumber\\
p(T,z)=\frac{\partial g}{\partial x}(X(T,z))\mathbb{E}[\delta_Z(z)|\mathcal{F}_T],
\end{array}
    \right.
\end{equation}
where
\begin{align}
\mu(t,z)&=-\frac{\partial H}{\partial x}(t,X(t,z),Y(t,z),u(t,z),p(t,z),q(t,z),r(t,z,.))\nonumber\\
&-\frac{\partial H}{\partial y}(t+\delta,X(t+\delta,z),Y(t+\delta,z),u(t+\delta,z),p(t+\delta,z),q(t+\delta,z),r(t+\delta,z,.))\mathbf{1}_{[0,T-\delta]}(t).
\end{align}
Let us now introduce a brief recall for the existence and uniqueness of the solution of a time-advanced BSDE. For more details see \cite{OSZ}.\\
Given a positive constant $\delta$, denote by $D([0,\delta],\mathbb{R})$ the space of all c\`adl\`ag paths from $[0,\delta]$ into $\mathbb{R}$.
For a path $X(.):\mathbb{R}_{+}\rightarrow\mathbb{R}$, $X_t$ will denote the function defined by
$X_t(s)=X(t+s)$ for $s\in[0,\delta]$.
Set $\mathcal{H}=L^2(\nu)$.\\
Let $F:\mathbb{R}_{+}\times\mathbb{R}\times\mathbb{R}
\times\mathbb{R}\times\mathbb{R}\times\mathcal{H}\times\mathcal{H}\times\Omega \rightarrow \mathbb{R}$ be a predictable function. Introduce the following Lipschitz condition:
There exists a constant $C$ such that
\begin{align}\label{lipsh}
&|F(t,p_1,p_2,q_1,q_2,r_1,r_2,\omega)-F(t,\bar{p}_1,\bar{p}_2,\bar{q}_1,\bar{q}_2,\bar{r}_1,\bar{r}_2,\omega)|\nonumber\\ &\leq C(|p_1-\bar{p}_1|+|p_2-\bar{p}_2|+|q_1-\bar{q}_1|+|q_2-\bar{q}_2|
+|r_1-\bar{r}_1|_{L^2{(\nu)}}+|r_2-\bar{r}_2|_{L^2{(\nu)}}),
\end{align}
where
\begin{equation}
|r|_{L^2{(\nu})}=(\int_{\mathbb{R}}r^2(\zeta)\nu(d\zeta))^{\frac{1}{2}}.
\end{equation}

Consider the following time-advanced BSDE in the unknown $\mathcal{F}_t$ adapted process $(p(t),q(t),r(t,\zeta))$
\begin{equation}\label{advbs}
\begin{cases}
dp(t)=E[F(t,p(t),p(t+\delta)1_{[0,T-\delta]}, q(t),q(t+\delta)1_{[0,T-\delta]}, r(t),r(t+\delta)1_{[0,T-\delta]})|\mathcal{F}_t]dt\\
+q(t)dB(t)+\int_{\mathbb{R}}r(t,\zeta)\tilde{N}(dt,d\zeta), t\in[0,T]\\
p(T)=G,
\end{cases}
\end{equation}
where $G$ is a given $\mathcal{F}_T$-measurable random variable such that $E[G^2]<\infty$.

Note that the time-advanced BSDE for the adjoint processes of the Hamiltonian is of this form.
For this type of time advanced BSDEs, we have the following result
\begin{theorem}\cite{OSZ}
Assume that condition \eqref{lipsh} is satisfied. Then the BSDE \eqref{advbs} has a unique solution $(p(t),q(t),r(t,\zeta))$ such that
\begin{equation}
E[\int_0^T \{ p^2(t)+q^2(t)+\int_{\mathbb{R}}r^2(t,\zeta)\nu(d\zeta)\}dt]<\infty.
\end{equation}
Moreover, the solution can be found by inductively solving a sequence of BSDEs backwards as follows.\\
Step 0. In the interval $[T-\delta, T ]$ we let $p(t), q(t),$ and $r(t, \zeta)$ be defined as the solution of the classical BSDE
\begin{equation}
\begin{cases}
dp(t) = F (t, p(t), 0, q(t), 0, r(t, \zeta), 0) dt + q(t) dB(t)+
\int_{\mathbb{R}}r(t, \zeta) \tilde{N}(dt, d\zeta), t\in [T-\delta, T ],\\
p(T ) = G.
\end{cases}
\end{equation}
Step $k, k\geq 1$. If the values of $(p(t), q(t), r(t,\zeta))$ have been found for $t\in[T-k\delta, T-(k-1)\delta]$ then, if $t\in[T-(k + 1)\delta, T-k\delta]$, the values of $p(t + \delta), q(t+\delta)$ and $r(t+\delta, \zeta)$
are known and, hence, the BSDE
\begin{equation}
\begin{cases}
dp(t) = E[F (t, p(t), p(t+\delta) , q(t), q(t + \delta), r(t), r(t+\delta))| \mathcal{F}_t ] dt\\
+ q(t) dB(t) +\int_{\mathbb{R}}r(t, \zeta) \tilde{N}(dt, d\zeta); t\in[T-(k +1)\delta, T-k\delta],\\
p(T-k\delta) = \text{ the value found in step $k-1$},
\end{cases}
\end{equation}
has a unique solution in $[T-(k + 1)\delta, T-k\delta]$.
We proceed like this until $k$ is such that $T-(k + 1)\delta\leq 0 < T-k\delta$ and then we solve the
corresponding BSDE on the interval $[0, T-k\delta]$.
\end{theorem}

Let us now give some conditions on $b$, $\sigma$ and $\gamma$ in order to ensure that a unique solution of  the advanced BSDE of the adjoint processs \eqref{advbs} exists.
In this case we have:
\begin{align}
&F(t,p_1,p_2,q_1,q_2,r_1,r_2,\omega)=-\frac{\partial H}{\partial x}(t,X(t,z),Y(t,z),u(t,z),p_1,q_1,r_1)\nonumber\\
&-\frac{\partial H}{\partial y}(t+\delta,X(t+\delta,z),Y(t+\delta,z),u(t+\delta,z),p_2,q_2,r_2)\mathbf{1}_{[0,T-\delta]}(t)\nonumber\\
&=-\mathbb{E}[\delta_Z(z)|\mathcal{F}_t] \frac{\partial f}{\partial x}(t,X(t,z),u(t,z),z)-\frac{\partial b}{\partial x}(t,X(t,z),Y(t,z),u(t,z),z)p_1\nonumber\\
& - \frac{\partial \sigma}{\partial x}(t,X(t,z),Y(t,z),u(t,z),z)q_1-\int_{\mathbb{R}}\frac{\partial\gamma}{\partial x}(t,X(t,z),Y(t,z),u(t,z),z,\zeta)r_1(\zeta)\nu(d\zeta)\nonumber\\
&-\{\frac{\partial b}{\partial y}(t,X(t+\delta,z),Y(t+\delta,z),u(t+\delta,z),z)p_2 + \frac{\partial \sigma}{\partial x}(t+\delta,X(t+\delta,z),Y(t+\delta,z),u(t+\delta,z),z)q_2\nonumber\\
&+\int_{\mathbb{R}}\frac{\partial\gamma}{\partial x}(t+\delta,X(t+\delta,z),Y(t+\delta,z),u(t+\delta,z),z,\zeta)r_2(\zeta)\nu(d\zeta)\}1_{[0,T-\delta]}(t).
\end{align}
Now let us verify condition \eqref{lipsh}:
\begin{align}
&|F(t,p_1,p_2,q_1,q_2,r_1,r_2,\omega)-F(t,\bar{p}_1,\bar{p}_2,\bar{q}_1,\bar{q}_2,\bar{r}_1,\bar{r}_2,\omega)|\nonumber\\
&\leq
|\frac{\partial b}{\partial x}(t,X(t,z),Y(t,z),u(t,z),z)||p_1-\bar{p}_1|\nonumber\\
&+|\frac{\partial \sigma}{\partial x}(t,X(t,z),Y(t,z),u(t,z),z)||q_1-\bar{q}_1|\nonumber\\
&+\int_{\mathbb{R}}|\frac{\partial \gamma}{\partial x}(t,X(t,z),Y(t,z),u(t,z),z)||r_1(\zeta)-\bar{r}_1(\zeta)|\nu(d\zeta)\nonumber\\
&+|\frac{\partial b}{\partial y}(t,X(t,z),Y(t,z),u(t,z),z)||p_2-\bar{p}_2|\nonumber\\
&+|\frac{\partial \sigma}{\partial y}(t,X(t,z),Y(t,z),u(t,z),z)||q_2-\bar{q}_2|\nonumber\\
&+\int_{\mathbb{R}}|\frac{\partial \gamma}{\partial y}(t,X(t,z),Y(t,z),u(t,z),z)||r_2(\zeta)-\bar{r}_2(\zeta)|\nu(d\zeta).
\end{align}

To guarantee that the Lipschitz condition is verified we assume that the derivatives of $b$, $\sigma$ and $\gamma$ with respect to $x$ and $y$ are bounded.

We can now state the first maximum principle for our problem \eqref{problem2}:

\begin{theorem}{[Sufficient-type maximum principle]}\\
Let $\hat{u} \in \mathcal{A}$, and denote the associated solution  of \eqref{eq3.3} and \eqref{eq4.2a} by $\hat{X}(t,z)$ and \\
$(\hat{p}(t,z),\hat{q}(t,z),\hat{r}(t,z,\zeta))$, respectively.  Assume that the following hold:
\begin{enumerate}
 \item $ x \rightarrow g(x,z)$ is concave for all $z$,
 \item $(x,y,u)\rightarrow H(t,x,y,u,z,\widehat{p}(t,z),\widehat{q}(t,z),\hat{r}(t,z,\zeta))$ is concave for all $t,z,\zeta$,
 \item $\sup_{w\in\mathbb{U}}H\big(t,\hat{X}(t,z),\widehat{Y}(t,z),w,\widehat{p}(t,z),\widehat{q}(t,z),\hat{r}(t,z,\zeta)\big)$\\
      $=H\big(t,\hat{X}(t,z),\widehat{Y}(t,z),\widehat{u}(t,z),\widehat{p}(t,z),\widehat{q}(t,z),\hat{r}(t,z,\zeta)\big)$ for all $t,z,\zeta.$
 \end{enumerate}
Then $\widehat{u}(\cdot,z)$ is an optimal insider control for Problem 3.3.
\end{theorem}

\dproof  By considering an increasing sequence of stopping times $\tau_n$ converging to $T$, we may assume that all local integrals appearing in the computations below are martingales and hence have expectation 0. See \cite{OS2}. We omit the details.\\
Choose arbitrary $u(.,z)\in\mathcal{A}$, and let the corresponding solution of \eqref{eq3.3} and \eqref{eq4.2a}  be $X(t,z)$ and $(p(t,x,z)$, $q(t,x,z)$, $r(t,x,z,\zeta))$ respectively.
For simplicity of notation we write\\
$f(t)=f(t,X(t,z),u(t,z))$,
$\widehat{f}(t)=f(t,x,\widehat{X}(t,z),\widehat{u}(t,z))$
 and similarly with $b$, $\widehat{b}$, $\sigma$, $\widehat{\sigma}$ and so on.\\
 Moreover put
\begin{equation}
 \hat{H}(t,z)=H(t,\widehat{X}(t,z),\widehat{Y}(t,z),\widehat{u}(t,z),z,\widehat{p}(t,z),\widehat{q}(t,z),\widehat{r}(t,z,.))
 \end{equation}
 and
 \begin{equation}
 H(t,z)=H(t,X(t,z), Y(t,z),u(t,z),z,\widehat{p}(t,z),\widehat{q}(t,z),\widehat{r}(t,z,.)).
 \end{equation}

 In the following we write  $\widetilde{f}=f-\widehat{f}$, $\widetilde{b}=b-\widehat{b}$, $\widetilde{X}=X-\widehat{X}$.\\
 Consider
 \begin{equation*}
    j(u(.,z))-j(\widehat{u}(.,z))=I_1+ I_2,
 \end{equation*}
 where
 \begin{equation}\label{eq4.7}
    I_1=\mathbb{E}[\int_0^T\{f(t)-\widehat{f}(t)\}\mathbb{E}[\delta_Z(z)|\mathcal{F}_t]dt], \quad I_2=\mathbb{E}[(g(X(T,z),z)-g(\hat{X}(T,z),z))\mathbb{E}[\delta_Z(z)|\mathcal{F}_T]].
 \end{equation}
 By the definition of $H$  and the concavity of $H$, we have
  \begin{align}\label{eq4.8}
    I_1 &= \mathbb{E}[\int_0^T\{H(t,z)-\widehat{H}(t,z)-\widehat{p}(t,z)\widetilde{b}(t,z) - \widehat{q}(t,z)\widetilde{\sigma}(t,z)\nonumber\\
   & -\int_{\mathbb{R}}\hat{r}(t,z,\zeta)\tilde{\gamma}(t,z,\zeta)\nu(d\zeta)\}dt]\nonumber\\
   &\leq\mathbb{E}\Big[\int_0^T(\frac{\partial \hat{H}}{\partial x}(t,z)\tilde{X}(t,z)+\frac{\partial \hat{H}}{\partial y}(t,z)\tilde{Y}(t,z)\nonumber\\
   &+\frac{\partial \hat{H}}{\partial u}(t,z)\tilde{u}(t,z)-\widehat{p}(t,z)\widetilde{b}(t,z) - \widehat{q}(t,z)\widetilde{\sigma}(t,z)\nonumber\\
   &-\int_{\mathbb{R}}\hat{r}(t,z,\zeta)\tilde{\gamma}(t,z,\zeta)\nu(d\zeta))dt\Big].
  \end{align}

Since $g$ is concave with respect to $x$ we have
\begin{align}
&(g(X(T,z),z)-g(\hat{X}(T,z),z))\mathbb{E}[\delta_Z(z)|\mathcal{F}_T]\nonumber\\
&\leq \frac{\partial g}{\partial x}(\hat{X}(T,z),z)\mathbb{E}[\delta_Z(z)|\mathcal{F}_T](X(T,z)-\hat{X}(T,z)),
\end{align}
and hence
\begin{align}\label{eq4.11}
   I_2 &\leq\mathbb{E}[\frac{\partial g}{\partial x}(\widehat{X}(T,z))\mathbb{E}[\delta_Z(z)|\mathcal{F}_T]\tilde{X}(T,z)]=\mathbb{E}[\widehat{p}(T,z)\widetilde{X}(T,z)] \\ \nonumber
    &= \mathbb{E}[\int_0^T \widehat{p}(t,z) d\widetilde{X}(t,z)+\int_0^T\widetilde{X}(t,z)d\widehat{p}(t,z)+\int_0^Td[\hat{p} , \tilde{X}]_t] \\ \nonumber
    &= \mathbb{E}[\int_0^T \{\widehat{p}(t,z)\widetilde{b}(t,z) +\widetilde{X}(t,z)\mathbb{E}[\hat{\mu}(t,z)|\mathbb{F}_t] \nonumber\\
   & + \widetilde{\sigma}(t,z)\widehat{q}(t,z)+\int_{\mathbb{R}}\tilde{\gamma}(t,z,\zeta)\hat{r}(t,z,\zeta)\nu(d\zeta)\}dt].\nonumber
    \end{align}

Combining \eqref{eq4.8} and \eqref{eq4.11} we obtain using the fact $X(t)=\hat{X}(t)=\xi(t)$ for all $t\in [-\delta,0]$

\begin{align}\label{eq4.11a}
j(u(.,z))-j(\widehat{u}(.,z))&\leq \mathbb{E}\Big[\int_0^T(\frac{\partial \hat{H}}{\partial x}(t,z)\tilde{X}(t,z)+\frac{\partial \hat{H}}{\partial y}(t,z)\tilde{Y}(t,z)\nonumber\\
&+\frac{\partial \hat{H}}{\partial u}(t,z)\tilde{u}(t,z)+\hat{\mu}(t,z)\tilde{X}(t,z))dt\Big]\nonumber\\
&=\mathbb{E}\Big[\int_{\delta}^{T+\delta}\{\frac{\partial \hat{H}}{\partial x}(t-\delta,z)+\frac{\partial \hat{H}}{\partial y}(t,z)\mathbf{1}_{[0,T]}(t)+\hat{\mu}(t-\delta,z)\}\tilde{Y}(t,z)dt\nonumber\\
&+\int_0^T\frac{\partial \hat{H}}{\partial u}(t,z)\tilde{u}(t,z)dt\Big],
\end{align}
where
\begin{equation}
\hat{\mu}(t-\delta,z)=\frac{\partial \hat{H}}{\partial x}(t-\delta,z)+\frac{\partial \hat{H}}{\partial y}(t,z)\mathbf{1}_{[0,T]}(t),
\end{equation}

then
\begin{align}
j(u(.,z))-j(\widehat{u}(.,z))&\leq\mathbb{E}\Big[\int_0^T\frac{\partial \hat{H}}{\partial u}(t,z)\tilde{u}(t,z)dt\Big]\leq 0.\nonumber
\end{align}

The last inequality holds because of the maximum condition of $H$.
Hence $j(u)\leq j(\hat{u})$.
Since $u\in\mathcal{A}$ was arbitrary, this shows that $\hat{u}$ is optimal.

\fproof
\section{A necessary-type maximum principle}

In some cases the concavity conditions of Theorem 4.1 do not hold.  In such situations a corresponding necessary-type maximum principle can be useful. For this, instead of the concavity conditions we need the following assumptions about the set of admissible control values:\\
\begin{itemize}
\item
$A_1$. For all $t_0\in [0,T]$ and all bounded $\mathcal{F}_{t_0}$-measurable random variables $\alpha(z,\omega)$, the control
$\theta(t,z, \omega) := \mathbf{1}_{[t_0,T ]}(t)\alpha(z,\omega)$ belongs to $\mathcal{A}$.
\item

For all $u; \beta_0 \in\mathcal{A}$ with $\beta_0(t,z) \leq K < \infty$ for all $t,z$  define
\begin{equation}\label{delta}
    \delta(t,z)=\frac{1}{2K}dist((u(t,z),\partial\mathbb{U})\wedge1 > 0,
\end{equation}
and put
\begin{equation}\label{beta(t,z)}
    \beta(t,z)=\delta(t,z)\beta_0(t,z).
\end{equation}
Then the control
\begin{equation*}
    \widetilde{u}(t,z)=u(t,z) + a\beta(t,z) ; \quad t \in [0,T]
\end{equation*}
belongs to $\mathcal{A}$ for all $a \in (-1, 1)$.\\
\item
$A3$. For all $\beta$ as in (\ref{beta(t,z)}) the derivative process
\begin{equation*}
    \chi(t,z):=\frac{d}{da}X^{u+a\beta}(t,z)|_{a=0}
\end{equation*}
exists, and belongs to $\mathbf{L}^2(\lambda\times \mathbf{P})$ and
\begin{equation}\label{d chi}
    \left\{
\begin{array}{l}
    d\chi(t,z) = [\frac{\partial b}{\partial x}(t,z)\chi(t,z)+\frac{\partial b}{\partial y}(t,z)\chi(t-\delta,z)+\frac{\partial b}{\partial u}(t,z)\beta(t,z)]dt\\
    +[\frac{\partial \sigma}{\partial x}(t,z)\chi(t,z)+\frac{\partial \sigma}{\partial y}(t,z)\chi(t-\delta,z)+\frac{\partial\sigma}{\partial u}(t,z)\beta(t,z)]d B(t) \\
    +\int_{\mathbb{R}}[\frac{\partial \gamma}{\partial x}(t,z,\zeta)\chi(t,z)+\frac{\partial \gamma}{\partial y}(t,z,\zeta)\chi(t-\delta,z)+\frac{\partial \gamma}{\partial u}(t,z,\zeta)\beta(t,z)]\tilde{N}(dt,d\zeta),\\
    \chi(t,z)  = 0 \quad  \forall t\in[-\delta,0].
   \end{array}
    \right.
\end{equation}
\end{itemize}

\begin{theorem}{[Necessary maximum principle]} \\
Let $\hat{u} \in \mathcal{A}$. Then the following are equivalent:
\begin{enumerate}
\item $\frac{d}{da}J((\hat{u}+a\beta)(.,z))|_{a=0}=0$ for all bounded $\beta \in \mathcal{A}$ of the form (\ref{beta(t,z)}).
\item $\frac{\partial H}{\partial u}(t,z)_{u=\hat{u}}=0$ for all $t\in[0,T].$
\end{enumerate}
\end{theorem}
\dproof
For simplicity of notation we write $u$ instead of $\hat{u}$ in the following. \\
By considering an increasing sequence of stopping times $\tau_n$ converging to $T$, we may assume that all local integrals appearing in the computations below are martingales and have expectation 0. See \cite{OS2}. We omit the details.\\
We can write
$$\frac{d}{da}J((u+a\beta)(.,z))|_{a=0}=I_1+I_2,$$\\
where
$$I_1=\frac{d}{da}\mathbb{E}[\int_0^Tf(t,X^{u+a\beta}(t,z),u(t,z)+a\beta(t,z),z)\mathbb{E}[\delta_Z(z)|\mathcal{F}_t]dt]|_{a=0}$$\\
and
$$I_2=\frac{d}{da}\mathbb{E}[g(X^{u+a\beta}(T,z),z)\mathbb{E}[\delta_Z(z)|\mathcal{F}_T]]|_{a=0}.$$
By our assumptions on $f$ and $g$ and by \eqref{eq4.1} we have
\begin{align}\label{iii1}
    I_1&=\mathbb{E}[\int_0^T\{\frac{\partial f}{\partial x}(t,z)\chi(t,z)+\frac{\partial f}{\partial y}(t,z)\chi(t-\delta,z)+\frac{\partial f}{\partial u}(t,z)\beta(t,z)\}\mathbb{E}[\delta_Z(z)|\mathcal{F}_t]dt]\nonumber\\
    &=\mathbb{E}[\int_0^T\{\frac{\partial H}{\partial x}(t,z)-\frac{\partial b}{\partial x}(t,z)p(t,z)-\frac{\partial \sigma}{\partial x}(t,z)q(t,z)-
    \int_{\mathbb{R}}\frac{\partial \gamma}{\partial x}(t,z,\zeta)r(t,z,\zeta)\nu(d\zeta)\}\chi(t,z)dt\nonumber\\
    &+\int_0^T\{\frac{\partial H}{\partial y}(t,z)-\frac{\partial b}{\partial y}(t,z)p(t,z)-\frac{\partial \sigma}{\partial y}(t,z)q(t,z)-
    \int_{\mathbb{R}}\frac{\partial \gamma}{\partial y}(t,z,\zeta)r(t,z,\zeta)\nu(d\zeta)\}\chi(t-\delta,z)dt\nonumber\\
    &+\int_0^T\frac{\partial f}{\partial u}(t,z)\beta(t,z)\mathbb{E}[\delta_Z(z)|\mathcal{F}_t]dt]
\end{align} 
and
\begin{equation}\label{iii2}
    I_2=\mathbb{E}[g'(X(T,z),z)\chi(T,z)\mathbb{E}[\delta_Z(z)|\mathcal{F}_T]]=\mathbb{E}[p(T,z)\chi(T,z)].
\end{equation}

By the It\^{o} formula
\begin{eqnarray}\label{iii22}
   I_2&=& \mathbb{E}[p(T,z)\chi(T,z)]=\mathbb{E}[\int_0^Tp(t,z)d\chi(t,z)+\int_0^T\chi(t,z)dp(t,z)+\int_0^Td[\chi,p](t,z)] \\ \nonumber
   &=& \mathbb{E}[\int_0^Tp(t,z)\{\frac{\partial b}{\partial x}(t,z)\chi(t,z)+\frac{\partial b}{\partial y}(t,z)\chi(t-\delta,z)+\frac{\partial b}{\partial u}(t,z)\beta(t,z)\}dt\\ \nonumber
   &-&\int_0^T\chi(t,z)\mathbb{E}[\mu(t,z)|\mathcal{F}_t]dt \\ \nonumber
     &+&\int_0^Tq(t,z) \{\frac{\partial \sigma}{\partial x}(t,z)\chi(t,z)+\frac{\partial \sigma}{\partial y}(t,z)\chi(t-\delta,z)+\frac{\partial\sigma}{\partial u}(t,z)\beta(t,z)\}dt \\ \nonumber
   &+&\int_0^T\int_{\mathbb{R}}\{\frac{\partial \gamma}{\partial x}(t,z,\zeta)\chi(t,z)+\frac{\partial \gamma}{\partial y}(t,z,\zeta)\chi(t,z)+\frac{\partial \gamma}{\partial u}(t,z,\zeta)\beta(t,z)\}r(t,z,\zeta)\nu(\zeta)dt].
  \end{eqnarray}

Summing (\ref{iii1}) and (\ref{iii22}) we get
\begin{align}
    &\frac{d}{da}J((u+a\beta)(.,y))|_{a=0}=I_1+I_2=\mathbb{E}[\int_0^T \chi(t,z)\{\frac{\partial H}{\partial x}(t,z)+\mu(t,z) \}dt\nonumber\\
    &+\int_0^T\chi(t-\delta,z)\frac{\partial H}{\partial y}(t,z)
    +\frac{\partial H}{\partial u}(t,z)\beta(t,z)dt]\nonumber\\
    &=\mathbb{E}[\int_0^T \chi(t,z)\{\frac{\partial H}{\partial x}(t,z)-\frac{\partial H}{\partial x}(t,z)-\frac{\partial H}{\partial y}(t+\delta,z)\mathbf{1}_{[0,T-\delta]}(t) \}dt\nonumber\\
    &+\int_0^T\chi(t-\delta,z)\frac{\partial H}{\partial y}(t,z)
    +\frac{\partial H}{\partial u}(t,z)\beta(t,z)dt]\nonumber\\
    &=\mathbb{E}[\int_0^T\frac{\partial H}{\partial u}(t,z)\beta(t,z)dt.
\end{align}

We conclude that
\begin{equation}
    \frac{d}{da}J((u+a\beta)(.,z))|_{a=0}=0,
\end{equation}
if and only if $\mathbb{E}[\int_0^T\frac{\partial H}{\partial u}(t,z)\beta(t,z)dt]=0$ for all bounded $\beta\in\mathcal{A}$ of the form (\ref{beta(t,z)}).\\
\noindent In particular, applying this to $\beta(t,z) = \alpha(z,\omega)\mathbf{1}_{[s,T]}(t)$ where $\alpha(z,\omega)$ is bounded and $\mathcal{F}_{t_0}$ measurable, $s\geq t_0$ we obtain
\begin{equation}
\mathbb{E}[\int_s^T\frac{\partial H}{\partial u}(t,z)dt\alpha ]=0.
\end{equation}
Differentiating with respect to $s$, we get
\begin{equation*}
   \mathbb{E}[ \frac{\partial H}{\partial u}(s,z)\alpha ]=0.
\end{equation*}
Since this holds  for all $s\geq t_0$ and for all $\alpha$, we conclude that
\begin{equation}
   \mathbb{E}[ \frac{\partial H}{\partial u}(s,z)|\mathcal{F}_{t_0}]=0.
\end{equation}
Assume that
\begin{equation}
s\mapsto\frac{\partial H}{\partial u}(s) \text{ is continuous}.
\end{equation}
By letting $s$ goes to $t_0$ then we get
\begin{equation}
E[\frac{\partial H}{\partial u}(t_0)|\mathcal{F}_{t_0}]=0.
\end{equation}
Then we deduce that $\frac{\partial H}{\partial u}(t_0)=0$ for all $t_0\in[0,T]$, hence for all $t\in[0,T], \frac{\partial H}{\partial u}(t)=0$.
\fproof

\section {Optimal inside harvesting in a population modelled by a delay equation}

Let us consider a population growth example.
We denote by $X(t,Z)$ a single population at time $t$ developing by a constant birth rate $\beta>0$ and a constant death rate $\alpha>0$ per inhabitant. $Z$ here is an inside information about the future environment for example coming from global warming.
In this model we take off immediately the dead from the population.
We denote by the constant $r>0$ the development period of each person ($r=9$ months for example).
A migration movement happens in this population and we assume that the global rate of the migration is distributed as a white noise $\sigma\dot{B}$.
We denote by $u(t,Z)$ the harvesting rate of the population.
The population change is given by the following SDDE:
\begin{equation}\label{eq0.5}
\begin{cases}
dX(t,Z)=(-\alpha X(t,Z)+\beta X(t-r,Z)-u(t,Z))dt+\sigma dB(t),\quad t\in[0,T]\\
X(s)=\eta(s), \quad -r\leq s\leq 0.
\end{cases}
\end{equation}
Here $\eta$ is a deterministic function.\\
We denote by $J$ the performance functional given by
\begin{equation}\label{eqperfo}
J(u)=E[\int_0^T e^{-\rho t}\frac{1}{\gamma}u^{\gamma}(t,Z)dt+\theta X(T,Z)],
\end{equation}
where $\theta$ is an $\mathcal{F}_T$-measurable strictly positive bounded random variable, $\rho>0$ and $\gamma\in(0,1)$.\\
Transforming the delayed stochastic control problem \eqref{eq0.5}-\eqref{eqperfo} into a $z$-parameterized $\mathbb{F}$-adapted delayed stochastic control problem  we get the following $z$-parameterized SDDE:
\begin{equation}\label{examtrans}
\begin{cases}
dX(t,z)=(-\alpha X(t,z)+\beta X(t-r,z)-u(t,z))dt+\sigma dB(t),\quad t\in[0,T]\\
X(s)=\eta(s), \quad -r\leq s\leq 0.
\end{cases}
\end{equation}
Let $\mathcal{A}$ be the set of admissible controls, we require that $u(t,z) > 0$, that $E[\int_0^Tu^2(t)dt]<\infty$. \\
Hypotheses $(E_1)$ in \cite{MS} is easily verified in this case:
since \eqref{examtrans} is a linear SDDE with constant coefficient $\alpha$ and $\beta$ then the Lipschitz condition is verified.
Also we have
$Y(t,z)=X(t-r,z)$ is $\mathcal{F}_t$-measurable then the drift function is $\mathcal{F}_t$-measurable.
Therefore we get form Theorem I.1 in \cite{MS} the existence and uniqueness of the solution of the parameterized SDDE \eqref{examtrans} such that $X(t, z)$
in $L^2(\Omega\times[0, T ])$ for each $z$.
 
The  transformed performance functional $J$ is given by
\begin{equation}\label{eq0.2}
J(u)=E[\int_0^T e^{-\rho t}\frac{1}{\gamma}u^{\gamma}(t,z)E[\delta_Z(z)|\mathcal{F}_t]dt+\theta X(T,z)E[\delta_Z(z)|\mathcal{F}_T]].
\end{equation}

In this case the Hamiltonian is given by
\begin{equation}
H(t,x,y,u,z,p,q)=E[\delta_Z(z)|\mathcal{F}_t]e^{-\rho t}\frac{1}{\gamma}u^{\gamma}(t,z)+(-\alpha x+\beta y-u)p+\sigma q.
\end{equation}
We have
\begin{equation}
\frac{\partial H}{\partial x}(t,x,y,u,z,p,q)= -\alpha p, \quad \frac{\partial H}{\partial y}(t,x,y,u,z,p,q)=\beta p,
\end{equation}
it follows that
\begin{equation}
\mu(t,z)=\alpha p(t,z)-\beta p(t+r,z)1_{[0,T-r]}(t).
\end{equation}
Therefore the advanced BSDE verified by the adjoint processes $(p(t),q(t))$ is given by
\begin{equation}\label{globsde}
\begin{cases}
dp(t,z)=(\alpha p(t,z)-\beta E[p(t+r,z)1_{[0,T-r]}|\mathcal{F}_t])dt + q(t,z)dB(t),\quad t\in[0,T]\\
p(T,z)=\theta E[\delta_Z(z)|\mathcal{F}_T].
\end{cases}
\end{equation}
Assume that $E[E[\delta_Z(z)|\mathcal{F}_T]^2]<\infty$.
Equation \eqref{globsde} is a linear advanced BSDE then it is easy to verify that the solution exists and it is unique.
We solve this BSDE recursively:\\
Step 1: If $t\in[T-r,T]$, the BSDE gets the form:
\begin{equation}
\begin{cases}
dp(t,z)=\alpha p(t,z)dt+q(t,z)dB(t), \quad T-r\leq t\leq T\\
p(T,z)=\theta E[\delta_Z(z)|\mathcal{F}_T].
\end{cases}
\end{equation}
This is a linear BSDE where the solution is given by
\begin{equation}\label{eq6.10sol}
p(t,z)=e^{\alpha (t-T)}E[\theta E[\delta_Z(z)|\mathcal{F}_T]|\mathcal{F}_t ], \quad t\in [T-r,T]
\end{equation}
with corresponding q(t,z) given by the martingale representation theorem.\\
Note that this solution is strictly positive since $\theta$ is a strictly positive bounded random
variable and $E[\delta_Z (z)|\mathcal{F}_T ]$ is strictly positive. Then $p(t, z)$ is strictly positive for all $t\in [T-r,T]$.\\
Step 2: If $t\in[T-2r,T-r]$ and $T-2r>0$, then  we get the BSDE:
\begin{equation}
\begin{cases}
dp(t,z)=(\alpha p(t,z)-\beta E[p(t+r,z)|\mathcal{F}_t])dt+q(t,z)dB(t)\\
p(T-r,z) \text{  is known from step 1}.
\end{cases}
\end{equation}
We have also $p(t+r,z)$ is known from step 1.
So this is a simple BSDE which can be solved for $p(t,z)$ and $q(t,z)$ for $t\in[T-2r, T-r]$.\\
The solution of this BSDE is given by
\begin{equation}\label{eq6.12 sol}
p(t,z)=e^{\alpha t}E[p(T-r)e^{-\alpha(T-r)}+\beta\int_t^{T-r}e^{-\alpha s}E[p(s+r,z)|\mathcal{F}_s]ds|\mathcal{F}_t],\quad t\in[T-2r, T-r]
\end{equation}
Note that this solution is strictly positive for all $t\in[T-2r, T-r]$.\\
We continue like this by induction until a step $j$ where $j$ satisfies $T-jr\leq 0<T-(j-1)r$.
With this method we end up with a solution $p(t,z)$ of the BSDE \eqref{globsde}.\\
The Hamiltonian $H$ can have a finite maximum over all $u$ only if
\begin{equation}
\frac{\partial H}{\partial u}(t)=E[\delta_Z(z)|\mathcal{F}_t]e^{-\rho t}u^{\gamma -1}(t,z)-p(t,z)=0.
\end{equation}
Then
\begin{equation}\label{u hat exp}
\hat{u}(t,z)=\frac{e^{\frac{\rho t}{\gamma-1}}}{(E[\delta_Z(z)|\mathcal{F}_t])^{\frac{1}{\gamma-1}}}(\hat{p}(t,z))^{\frac{1}{\gamma-1}}.
\end{equation}

Let us now verify the admissibility condition of $\hat{u}(t,z)$.
\begin{itemize}
\item The solution $p(t,z)$  of equation \eqref{globsde} is strictly positive since $\theta$ is a strictly positive bounded random
variable and $E[\delta_Z(z)|\mathcal{F}_t]>0$.
Then $\hat{u}(t,z)$ in \eqref{u hat exp} is strictly positive.
\item 
We will prove that $E[\int_{0}^T\hat{u}^2(t)dt]<\infty$ by steps.\\
Since $\theta$ is bounded positive random variable being away from 0 i.e there exist $\theta_1$, $\theta_2\in{\mathbb{R}_+}$ such that 
\begin{equation}
\theta_1\leq \theta(\omega)\leq \theta_2 ,
\end{equation}
therefore we have
\begin{itemize}
\item i) For $t \in [T-r,T]$ we have by \eqref{eq6.10sol} 
\begin{equation}
p(t,z) \geq e^{\alpha (t-T)}\theta_1E[\delta_{Z}(z) | \mathcal{F}_t].
\end{equation}

\item ii) Then for $t \in [T-2r,T-r]$ we have by \eqref{eq6.12 sol} and (i)
\begin{equation}
p(t,z) \geq e^{-\alpha(T-r-t)} E[p(T-r) | \mathcal{F}_t] \geq e^{-\alpha(T-t)}\theta_1 E[ E[\delta_{Z}(z) | \mathcal{F}_{T-r} | \mathcal{F}_t]] \geq e^{-\alpha(T-t)}\theta_1 E[\delta_{Z}(z) | \mathcal{F}_t].
\end{equation}

\item iii) Proceeding like this by induction we end up with
\begin{equation}
p(t.z) \geq e^{-\alpha(T-t)}\theta_1 E[\delta_{Z}(z) | \mathcal{F}_t], \forall t \in [0,T].
\end{equation}

\end{itemize}

We conclude that 
\begin{equation}
E[\int_0^T u^2(t) dt] \leq e^{-\frac{2\alpha(T-t)}{\gamma-1}}\theta_1^{\frac{2}{\gamma-1}} E[ \int_0^T E[\delta_{Z}(z) | \mathcal{F}_t]^{2/(1-\gamma)} E[\delta_{Z}(z) | \mathcal{F}_t]^{2/(\gamma-1)}dt =e^{-\frac{2\alpha(T-t)}{\gamma-1}}\theta_1^{\frac{2}{\gamma-1}} T.
\end{equation}

\end{itemize}
\begin{theorem}
Assume that $\gamma\in(0,1)$ and $E[E[\delta_Z(z)|\mathcal{F}_T]^2]<\infty$.
The harvesting control in the population growth example  given by \eqref{eq0.5}-\eqref{eqperfo} is given by 
\begin{equation}
\hat{u}(t,Z)=\frac{e^{\frac{\rho t}{\gamma-1}}}{(E[\delta_Z(z)|\mathcal{F}_t]_{z=Z})^{\frac{1}{\gamma-1}}}(\hat{p}(t,Z))^{\frac{1}{\gamma-1}},
\end{equation}
where $p(t,z)$ is the solution of the advanced BSDE \eqref{globsde}.
\end{theorem}

\section {Optimal insider portfolio in a financial market with delay}

Consider the following SDDE:
\begin{equation}\label{BI1}
\begin{cases}
dX(t,Z)=X(t,Z)\pi(t,Z)[b(t)\theta(t,Z)dt+\sigma(t)\theta(t,Z)dB(t)], t\in[0,T]\\
X(t)=\xi(t), t\in[-r,0].
\end{cases}
\end{equation}
 Here $b, \sigma$ and $\xi$ are deterministic bounded functions with $\xi(0)=1$ and  $\theta(t,Z)=\frac{X(t-r,Z)}{X(t,Z)}1_{t<\tau_0}$ where
\begin{equation}
\tau_0=\inf\{t>0, X(t,Z)=0 \}.
\end{equation}
The process $\pi$ is our control process and assumed to be $\mathbb{H}$ adapted.\\
Note that equation \eqref{BI1} is well defined since $B$ is a semimartingale under the filtration $\mathbb{H}$ as discussed in Section 2.
Since we are working via Hida-Malliavin calculus we will interpret it  using forward integral definition.
Here a brief recall of the definition of forward integral:
\begin{definition}
We say that a stochastic process $\phi = \pi(t), t\in [0, T ]$, is forward integrable
(in the weak sense) over the interval $[0, T ]$ with respect to $B$ if there exists a process $I =I(t), t \in [0, T ]$, such that
\begin{equation}
\sup_{t\in[0,T]} \int_0^t |\phi(s)\frac{B(s + \epsilon)-B(s)}{\epsilon}ds-I(t)|\rightarrow 0, \quad \epsilon\rightarrow 0^+
\end{equation}
in probability.\\

 In this case we write
$I(t) :=\int_0^t\phi(s)d^-B(s), t \in [0, T ]$,
and call $I(t)$ the forward integral of $\phi$ with respect to $B$ on $[0, t]$.
\end{definition}

For $\pi$ to be admissible in the set $\mathcal{A}$ of admissible $\mathbb{H}$ adapted controls, we require that
\begin{equation}\label{admissible requirement}
E[\int_0^{T\wedge\tau_0}\frac{X^2(s-r,Z)}{X^2(s,Z)} \pi^2(s,Z)ds]<\infty.
\end{equation}
Condition \eqref{admissible requirement} guarantees the existence and uniqueness of the solution of \eqref{BI1}.\\
 The solution of this SDDE is given by
 \begin{align}
 X(t,Z)&=\exp(\int_0^{t\wedge\tau_0}(b(s)\frac{X(s-r,Z)}{X(s,Z)}\pi(s,Z)-\frac{1}{2}\sigma^2(s)\frac{X^2(s-r,Z)}{X^2(s,Z)}\pi^2(s,Z))ds\nonumber\\
 &+\int_0^{t\wedge\tau_0}\sigma(s)\frac{X(s-r,Z)}{X(s,Z)}\pi(s,Z)dB(s)), t\in[0,T].
 \end{align}
 \textbf{Problem}
 Our goal is to find $\pi^*\in \mathcal{A}$ which maximizes the following expected utility logarithmic function
 \begin{equation}\label{eq7.5new}
 \sup_{\pi\in\mathcal{A}}E\Big[\ln(X^{\pi}(T\wedge\tau_0,Z))\Big]=E\Big[\ln(X^{\hat{\pi}}(T\wedge\tau_0,Z))\Big].
 \end{equation}
We have
 \begin{align}
&E\Big[\ln(X^{\pi}(T\wedge\tau_0,Z))\Big]=E\Big[\int_0^{T\wedge\tau_0}\Big(b(s)\frac{X(s-r,Z)}{X(s,Z)}\pi(s,Z)-\frac{1}{2}\sigma^2(s)\frac{X^2(s-r,Z)}{X^2(s,Z)}\pi^2(s,Z)\nonumber\\
&+\sigma(s)E[D_s(\frac{X(s-r,Z)}{X(s,Z)}\pi(s,Z))|\mathcal{F}_s]\Big)ds\Big]\nonumber\\
&=E\Big[\int_0^{T\wedge\tau_0}E\Big[\Big(b(s)\frac{X(s-r,Z)}{X(s,Z)}\pi(s,Z)-\frac{1}{2}\sigma^2(s)\frac{X^2(s-r,Z)}{X^2(s,Z)}\pi^2(s,Z)\nonumber\\
&+\sigma(s)D_s(\frac{X(s-r,Z)}{X(s,Z)}\pi(s,Z))\Big)|\mathcal{F}_s\Big]ds\Big].
\end{align}
Here and in the following we use the notation
\begin{equation}
D_s\pi(s) := D_{s^+}\pi(s) := \lim_{t\rightarrow s^+} D_t \pi(s)
\end{equation}
and the following result:(see \cite{DOP})\\

Let $\phi$ be a c\`agl\`ad and forward integrable process in $L^2([0,T]\times\Omega )$ then 
\begin{equation}
E[\int_0^T\phi(s)d^-B(s)]=E[\int_0^TE[D_{s^ +}\phi(s)|\mathcal{F}_s ]ds].
\end{equation}
Let
\begin{equation}
J(\pi)=E\Big[\Big(b(s)\frac{X(s-r,Z)}{X(s,Z)}\pi(s,Z)-\frac{1}{2}\sigma^2(s)\frac{X^2(s-r,Z)}{X^2(s,Z)}\pi^2(s,Z)+\sigma(s)D_s(\frac{X(s-r,Z)}{X(s,Z)}\pi(s,Z))\Big)|\mathcal{F}_s\Big].
\end{equation}
Applying the Donsker delta functional to the previous equation we get 
\begin{align}
J(\pi)&=\int_{\mathbb{R}}\Big(b(s)\frac{X(s-r,z)}{X(s,z)}\pi(s,z)E[\delta_Z(z)|\mathcal{F}_s]-\frac{1}{2}\sigma^2(s)\frac{X^2(s-r,z)}{X^2(s,z)}\pi^2(s,z)E[\delta_Z(z)|\mathcal{F}_s]\nonumber\\
&+\sigma(s)\frac{X(s-r,z)}{X(s,z)}\pi(s,z)E[D_s\delta_Z(z)|\mathcal{F}_s]\Big)dz.
\end{align}
We can maximize this over $\pi(s,z)$ for each $s,z$.
Then we get that
\begin{equation}
\hat{\pi}(s,z)=\frac{X(s,z)}{\sigma(s)X(s-r,z)}\frac{E[D_s\delta_Z(z)|\mathcal{F}_s]}{E[\delta_Z(z)|\mathcal{F}_s]}+\frac{b(s)}{\sigma^2(s)}\frac{X(s,z)}{X(s-r,z)}, s\in [0,T\wedge\tau_0].
\end{equation}
Therefore we get
\begin{align}
\hat{\pi}(s,Z)&=\frac{X(s,Z)}{\sigma(s)X(s-r,Z)}\frac{E[D_s\delta_Z(z)|\mathcal{F}_s]_{z=Z}}{E[\delta_Z(z)|\mathcal{F}_s]_{z=Z}}+\frac{b(s)}{\sigma^2(s)}\frac{X(s,Z)}{X(s-r,Z)}, s\in [0,T\wedge\tau_0]\nonumber\\
&=\frac{X(s,Z)}{\sigma(s)X(s-r,Z)}\Phi(s,Z)+\frac{b(s)}{\sigma^2(s)}\frac{X(s,Z)}{X(s-r,Z)}, s\in [0,T\wedge\tau_0],
\end{align}
where
\begin{equation}
\Phi(s,Z)=\frac{E[D_s\delta_Z(z)|\mathcal{F}_s]_{z=Z}}{E[\delta_Z(z)|\mathcal{F}_s]_{z=Z}}.
\end{equation}
Replacing the expression of $\hat{\pi}$ in $\ln(X^{\pi}(T\wedge\tau_0,Z))$ we get
\begin{align}
 \ln(X^{\hat{\pi}}(T\wedge\tau_0,Z))&=\int_0^{T\wedge\tau_0}\Big(b(s)(\frac{\Phi(s,Z)}{\sigma(s)}+\frac{b(s)}{\sigma^2(s)})-\frac{1}{2}\sigma^2(s) (\frac{\Phi(s,Z)}{\sigma(s)}+\frac{b(s)}{\sigma^2(s)})^2  \Big)ds\nonumber\\
 &+\int_0^{T\wedge\tau_0}\sigma(s)(\frac{\Phi(s,Z)}{\sigma(s)}+\frac{b(s)}{\sigma^2(s)})dB(s)>-\infty,
\end{align}
and hence $X^{\hat{\pi}}(T\wedge\tau_0,Z)>0$ a.s.\\
This is only possible if $\tau_0>T$ a.s, which means that $X^{\hat{\pi}}(t)>0$ for all $t\in[0,T]$ a.s and our optimal indeed in the whole interval $[0,T]$.
Let us now verify the admissibility condition:
\begin{align}
E[\int_0^{T}\frac{X^2(s-r,Z)}{X^2(s,Z)} \hat{\pi}^2(s,Z)ds]=E[\int_0^{T}(\frac{E[D_s\delta_Z(z)|\mathcal{F}_s]_{z=Z}}{E[\delta_Z(z)|\mathcal{F}_s]_{z=Z}}+\frac{b(s)}{\sigma^2(s)})^2ds].
\end{align}
This previous quantity is finite if we suppose that
\begin{equation}
E[\int_0^{T}(\frac{E[D_s\delta_Z(z)|\mathcal{F}_s]_{z=Z}}{E[\delta_Z(z)|\mathcal{F}_s]_{z=Z}})^2ds]<\infty.
\end{equation}
Let us now resume the previous result in the following Theorem:
\begin{theorem}
Assume that \begin{equation}
E[\int_0^{T}(\frac{E[D_s\delta_Z(z)|\mathcal{F}_s]_{z=Z}}{E[\delta_Z(z)|\mathcal{F}_s]_{z=Z}})^2ds]<\infty.
\end{equation}
Then the optimal portfolio $\hat{\pi}$ with respect to the logarithmic utility for an insider in the delay market \eqref{BI1} and with inside information \eqref{eq1.1} is given by

\begin{align}
\hat{\pi}(s,Z)&=\frac{X(s,Z)}{\sigma(s)X(s-r,Z)}\frac{E[D_s\delta_Z(z)|\mathcal{F}_s]_{z=Z}}{E[\delta_Z(z)|\mathcal{F}_s]_{z=Z}}+\frac{b(s)}{\sigma^2(s)}\frac{X(s,Z)}{X(s-r,Z)}, s\in [0,T]
\end{align}
\end{theorem}
In the following Corollary we treat a particular case:
\begin{coro}
Suppose that $Z=B(T_0)$, where $T_0>T$. In this case we have
\begin{equation}
\frac{E[D_s\delta_Z(z)|\mathcal{F}_s]_{z=Z}}{E[\delta_Z(z)|\mathcal{F}_s]_{z=Z}}=\frac{B(T_0)-B(s)}{T_0-s},
\end{equation}
and
\begin{equation}
E[\int_0^T(\frac{B(T_0)-B(s)}{T_0-s})^2ds]=\int_0^T\frac{ds}{T_0-s}<\infty \text{  since  } T_0>T.
\end{equation}
Then
\begin{equation}
\hat{\pi}(s,B(T_0))=\frac{X(s,B(T_0))}{\sigma(s)X(s-r,B(T_0))}\frac{B(T_0)-B(s)}{T_0-s}+\frac{b(s)}{\sigma^2(s)}\frac{X(s,B(T_0))}{X(s-r,B(T_0))}, s\in [0,T].
\end{equation}
\end{coro}

\subsection{Viability of a market with delay}
In this Subsection we study the viability of the financial market \eqref{BI1}-\eqref{eq7.5new}. For this we define what is a viable market:

\begin{definition}
The martket \eqref{BI1}-\eqref{eq7.5new} is called viable if
\begin{equation}
\sup_{\pi\in\mathcal{A}}E[\ln X^{\pi}(T)]<\infty.
\end{equation}
\end{definition}

In the no delay case $(r=0)$ we have that
\begin{equation}
\hat{\pi}(s,B(T_0))=\frac{B(T_0)-B(s)}{\sigma(s)(T_0-s)}+\frac{b(s)}{\sigma^2(s)}, s\in [0,T].
\end{equation}
In this case we know from \cite{PK} that for $T_0=T, E[\ln(X^{\hat{\pi}}(T,B(T)))]$ is infinite. In the case of delay for the market \eqref{BI1}, one can show also that for $T_0=T$, we get $E[\ln(X^{\hat{\pi}}(T,B(T)))]=\infty$.
In fact we have
\begin{align}\label{eq7.20BI}
E[\ln(X^{\hat{\pi}}(T,B(T)))]&=E[\int_0^T\{b(s)(\frac{\Phi(s,Z)}{\sigma(s)}+\frac{b(s)}{\sigma^2(s)})-\frac{1}{2}\sigma^2(s)(\frac{\Phi(s,Z)}{\sigma(s)}+\frac{b(s)}{\sigma^2(s)})^2\}ds]\nonumber\\
&+E[\int_0^T\sigma(s)(\frac{\Phi(s,Z)}{\sigma(s)}+\frac{b(s)}{\sigma^2(s)})dB(s)]\nonumber\\
&=E[\int_0^T\{b(s)(\frac{\Phi(s,Z)}{\sigma(s)}+\frac{b(s)}{\sigma^2(s)})-\frac{1}{2}\sigma^2(s)(\frac{\Phi(s,Z)}{\sigma(s)}+\frac{b(s)}{\sigma^2(s)})^2\}ds]\nonumber\\
&+E[\int_0^TD_s\Phi(s,Z)ds],
\end{align}

where $b$ and $\sigma$ are deterministic bounded.
We have
\begin{equation}
D_s\Phi(s,Z)=D_s(\frac{B(T)-B(s)}{T-s})=\frac{1}{T-s}1_{[0,T]}.
\end{equation}
Then
\begin{align}
E[\ln(X^{\hat{\pi}})(T,B(T))]&=E[\int_0^T\{b(s)(\frac{B(T)-B(s)}{\sigma(s)(T-s)}+\frac{b(s)}{\sigma^2(s)})-\frac{1}{2}\sigma^2(s)(\frac{B(T)-B(s)}{\sigma(s)(T-s)}+\frac{b(s)}{\sigma^2(s)})^2\}ds]\nonumber\\
&+E[\int_0^T\frac{1}{T-s}1_{[0,T]} ds]\nonumber\\
&=E[\int_0^T\{ -\frac{1}{2}\frac{(B(T)-B(s))^2}{(T-s)^2}-b(s)\frac{B(T)-B(s)}{\sigma(s)(T-s)}+\frac{1}{2}\frac{b^2(s)}{\sigma^2(s)}\}ds]+E[\int_0^T\frac{1}{T-s} ds]\nonumber\\
&=\frac{1}{2}\int_0^T\frac{ds}{(T-s)}+\int_0^T\frac{b^2(s)}{\sigma^2(s)}ds=\infty.
\end{align}
So we conclude that even in a market with delay, the market is not viable when the inside information $Z=B(T)$.

\end{document}